\documentclass[preprint,number,12pt]{elsarticle}
\usepackage{amsmath,amssymb,graphicx,longtable}

\usepackage{lineno}

\begin{document}

\begin{frontmatter}

\title{Separation of carbon dioxide from flue emissions using Endex principles}

\author[anu]{R. Ball\corref{cor1}}
\ead{Rowena.Ball@anu.edu.au}

\author[calix]{M. G. Sceats\corref{cor2}}
\ead{mgsceats@calix.com.au}

\cortext[cor1]{Principal corresponding author}
\cortext[cor2]{Corresponding author}

\address[anu]{Mathematical Sciences Institute and Research School of Physics and Engineering,\\ The Australian National University, Canberrra ACT 0200 Australia }

\address[calix]{Calix Limited, 828 Pacific Highway Gordon NSW 2072 Australia}
\date{\normalsize{10 February 2010}}

\begin{abstract}
In an Endex reactor endothermic and exothermic reactions are directly thermally cou­pled and kinetically matched to achieve intrinsic thermal stability, efficient conversion, autothermal operation, and minimal heat losses. Applied to the problem of in-line carbon dioxide separation from flue gas, Endex principles hold out the promise of effecting a CO$_2$-capture technology of unprecedented economic viability. In this work we describe an Endex Calcium Looping reactor, in which heat released by chemisorption of carbon dioxide onto calcium oxide is used directly to drive the reverse reaction, yielding a pure stream of CO$_2$ for compression and geosequestration. In this initial study we model the proposed reactor as a continuous-flow dynamical system in the well-stirred limit, compute the steady states and analyse their stability properties over the operating parameter space, flag potential design and opera­tional challenges, and suggest an optimum regime for effective operation.
\end{abstract}


\end{frontmatter}
\section{Introduction }

In this paper we introduce a method of reactive endothermic-exothermic coupling, known as an Endex system \cite{Gray:1999,Ball:1999}, for the in-line separation of carbon dioxide from flue and fuel gas emissions.  Endex  principles are applied to obtain a highly efficient modification of the Calcium Looping separation technique, in which calcium oxide (CaO, or lime) is used to scrub CO$_2$ from a flue gas \cite{Anthony:2009}. The Endex process involves direct thermal coupling of a carboniser reactor segment, in which the lime sorbent reacts with CO$_2$ to produce calcium carbonate (CaCO$_3$) and a scrubbed gas effluent stream, and a calciner reactor segment to which the loaded sorbent is transported and where the sorbent is regenerated and reinjected into the carboniser segment, with production of a pure CO$_2$ gas stream. The same methodology can be applied to other sorbents and to fuel gas mixtures such as syngas.  Since the Endex Calcium Looping reactor is a thermoreactive system an important first step in a project to build and operate a demonstration plant is to assess the thermal stability of the system. In this work we model the reactor as a coupled dynamical system in the well-stirred limit and map the linear stability of the steady states over the operational parameter space.

Technologies for scrubbing carbon dioxide from flue gas to reduce greenhouse gas emissions must  satisfy the following criteria, at least, in order to be accepted by society and implemented by industry: they must (1) effect CO$_2$ capture for less than $\sim$US\$30 per tonne with an increased cost of electricity of less than about 15\%, (2) preferably be retrofittable to existing fossil-fueled plants, (3) not cause additional environmental harm, and (4) operate safely. The best available technologies may not be able to comply with all of these requirements \cite{Yang:2008}. There is a drive to develop second generation technologies which have been purpose-designed to meet the stringent demands of CO$_2$ capture, rather than adapt existing technologies.  

Calcium Looping was first proposed by Heesink and  Temmink in 1994 \cite{Heesink:1994} for post-combustion CO$_2$ emissions reduction, although its use on an industrial scale for removal of CO$_2$ from syngases dates from the 1960s \cite{Curran:1967}. Shimizu \textit{et al.} \cite{Shimizu:1999} proposed that oxyfuel combustion of additional fossil fuel be used to drive  the high temperature calcination process. Abanades \textit{et al.}  \cite{Abanades:2007} showed that the CO$_2$--lime chemical bond strength  is such that the thermal energy required is about 30--40\% of the thermal energy of the plant that produces the flue gas. This additional consumption of energy results in higher CO$_2$ production, but this is offset by use of the heat liberated from the carboniser at about 873--1023\,K to produce additional power.  

The major drawback of Calcium Looping is the rapid loss of carbonation capacity  due to sintering of the lime sorbent. Sintering occurs at the high temperatures used in the conventional calcination cycle, about 1123--1223\,K in an atmosphere of carbon dioxide \cite{Borgwardt:1989a,Borgwardt:1989b} and, more importantly, reactive sintering results in irreversible loss of the mesopore surface area due to pore filling \citep{Gonzales:2008}. 

Sintering was recognised as a problem by Curran \textit{et al.} \cite{Curran:1967} and more recently by Silaban and  Harrison \cite{Silaban:1995} and has been the subject of intensive research since. For example, Abanades and Alvarez \cite{Abanades:2003} showed that the sorption capacity degrades from a capacity of about 80\% of the theoretical limit of 0.78 kg of CO$_2$ per kg of CaO to about 17\% after about 10 cycles, and reactivation studies were carried out by Fennell \textit{et al.} \cite{Fennell:2007} and Manovic and Anthony \cite{Manovic:2008}. Even though MacKenzie \textit{et al.} \cite{MacKenzie:2007} showed that CO$_2$ capture using Calcium Looping with sintering meets criteria (1)--(3) above, sorbent sintering remains the largest barrier to adoption of Calcium Looping. In the Endex Calcium Looping process described in the present paper the sintering problem is neatly sidestepped by inversion of the usual temperature difference between carboniser and calciner and by using light carbonation over fast cycling.

The Endex approach to thermal control and heat recovery was first described and analysed in mathematical terms by Gray and Ball \cite{Gray:1999,Ball:1999}. It involves using the heat generated by an exothermic reaction to drive an endothermic reaction directly, in real time. The reactions are thermally coupled via direct heat transfer or mass transfer or both, and for the reactions to effectively ``see'' each other they must be matched kinetically Since the kinetic parameters (activation energies and pre-exponential factors) of a selected Endex couple do not match in general, kinetic matching is achieved (to an acceptable approximation) by manipulation of residence times.  Several studies of specific Endex systems have been carried out since with promising results, although none has referred to the original works of Gray and Ball \cite{Gray:1999,Ball:1999}. Terms such as ``recuperative coupling'' and ``co-current thermally coupled reactor'' are used for Endex systems and principles in those works, which are reviewed in a recent study of a thermally coupled reactor \cite{Khademia:2009}.

Sceats \cite{Sceats:2009a} and Sweeney and Sceats \cite{Sweeney:2009} proposed that an Endex configuration could be used in which the temperature of the calciner is held below that of the carboniser by thermally coupling the reactor segments.  This configuration is counter-intuitive because it is usually assumed that the endothermic, bond-breaking, process in the calciner would occur at a higher temperature than that of the carboniser, and all previous approaches to Calcium Looping are based on this conventional configuration.  Inversion of the temperature difference between carboniser and calciner is achieved by control of the pressure.  Operation using this new configuration has three consequences: the lower temperature of the calciner reduces thermal sintering, partial carbonation reduces the loss of surface area by pore-filling, and the unreacted sorbent promotes heat transfer.  The intent was to greatly reduce the impact of thermal and reactive sintering, and remove this barrier to adoption of Calcium Looping.  The most desirable property of the Endex configuration is that CO$_2$ separation can be achieved, in principle, without the need for additional heat.  Sweeney and Sceats \cite{Sweeney:2009} conclude that the Endex configuration has the potential to significantly reduce the cost of CO$_2$ capture below that of current technologies, principally because of this property.

\begin{figure}[ht]
\centerline{
  \includegraphics[scale=0.5]{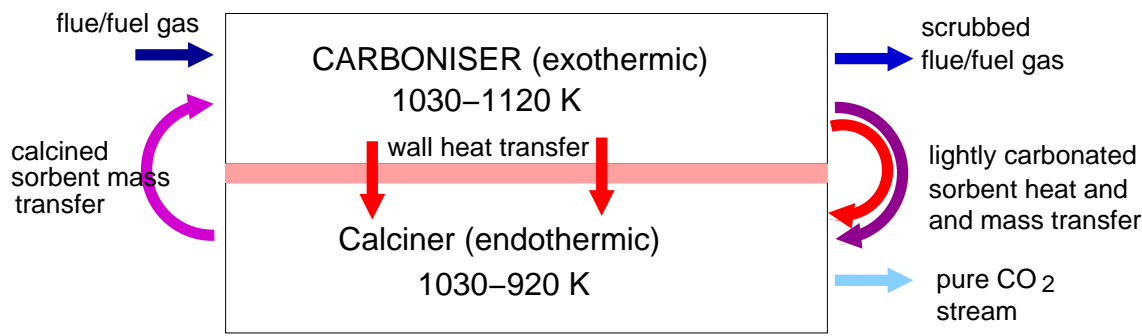}}
  \caption[]{\label{fig1} In the Endex configuration the calciner temperature is lower than the carboniser temperature. Heat transfer by sorbent mass flow is assisted by direct wall heat transfer, and the CO$_2$ pressure in the calciner is low. Light carbonation only is required; this and the lower calcination temperature sidestepps the sintering problem.}
\end{figure}

A schematic of the Endex Calcium Looping reactor is given in figure \ref{fig1}, where the wall and mass thermal couplings are indicated.

Thermal coupling of reactive systems creates additional  nonlinearities in the dynamics.  In the Endex reactor configuration for CaO Looping, CO$_2$ separation can be considered as a gas switch in which the CO$_2$ pressure in the calciner, for example, can be used to control the CO$_2$ output mass flow.  However, reactor coupling can amplify perturbations with adverse effects, and a more detailed analysis of the Endex configuration is required that deals with its response to perturbations, and the start-up and shut down processes. This analysis is performed in this paper.  In section \ref{sec2} we review the properties of Endex systems, and outline the rationale for stability analysis of thermoreactive systems generally and of the proposed Endex configuration in partic­ular. The methodology used is described in section \ref{sec3}, and the dynamical model is presented. In section \ref{sec4} the results are presented and discussed for: standalone carboniser, Endex carboniser/­calciner with sorbent cycling but without wall thermal coupling, Endex carboniser-calciner with both sorbent cycling and wall thermal coupling, and transient analysis (start-up and shutdown dynamics). We summarize these results in section \ref{sec5}.

\section{Carbonation, calcination, and Endex principles \label{sec2}} 

The chemisorption of carbon dioxide onto calcium oxide is highly exothermic, releasing  $170 $\,kJ/mol CO$_2$ at 1073\,K. Housed in a reactor where thermostatting is necessarily imperfect such a reaction may become thermally unstable and either self-quench or self-heat uncontrollably. Both situations are obviously undesirable from safety and economic considerations. 
 
Traditionally, reactors housing exothermic reactions are fitted with thermostatted cooling jackets, but the problem of maintaining thermal stability becomes much more difficult (and expensive) as the reactor is scaled up because the heat removal rate is linear in the temperature and scales as $L^2$ while the reactive heat generation rate is exponential in the temperature and scales as $L^3$, where $L$ is a characteristic reactor dimension.  This scale-up problem is particularly relevant to reactive flue gas CO$_2$ capture, because the emitter is typically a large fossil-fueled power~station. 

Reaction heat that is removed by the cooling system in traditional chemical reactor plants is either dissipated into the environment, or sometimes, in large and sophisticated plants, partially recovered by an indirect route and used for another purpose.

It is important to appreciate the potential for adiabatic, or thermally insulated, operation of an Endex configured reactor. A conventional chemical reactor is almost never run adiabatically (unless it is a bomb calorimeter) because for most common industrial reactions the adiabatic temperature rise for full conversion is dangerously high. Purpose-built cooling systems are usually necessary, which can become technically very elaborate, and expensive, when the reactor is large. The concept of an \textit{ideal adiabatic} Endex reactor, where the additional conservation condition of enthalpy flux conservation holds to a good approximation, suggests the very appealing possibility that scaling problems may be eliminated, in the same way that surface-to-volume ratios are irrelevant in the single conventional adiabatic reactor, while good conversion, thermal safety, and direct recovery of reaction heat are also achieved. 

An Endex system is a coupled nonlinear dynamical system, and therefore has the potential to exhibit thermal instabilities and is capable of more complex behavior than a single insulated or thermostatted exothermic reaction system. 
 However, in the works of Gray and Ball \cite{Gray:1999,Ball:1999} it was shown that a general Endex-configured reaction system can operate stably, autothermally, and economically, achieving almost full recovery of chemical energy and co-production of valuable products. The tradeoff is that the stability regime must be mapped for specific Endex systems in order to avoid or control (or even, perhaps, exploit) the coupled relaxation oscillator dynamics that are inherent to coupled dynamical systems. 

In the proposed Endex reactor  the exothermic carbonation reaction \linebreak CaO~+~CO$_2$~$\rightleftarrows$~CaCO$_3$ is thermally coupled and kinetically matched with the reverse  reaction, the endothermic calcination reaction \linebreak CaCO$_3$ $\rightleftarrows$ CaO + CO$_2$. To show that, in principle, this system can be operated effectively within a large margin of thermal safety we have carried out a linear stability analysis on the steady-state solutions of a dynamical model for the reactor.

Stability analysis is a valuable reactor system design tool, and an essential step in the design of thermoreactive systems. Typically we are interested in running a continuous-flow thermoreactive system at a particular steady state, or set-point. Stability analysis tells us whether small perturbations around the steady state---which are inevitable in any real system---will decay and settle back onto the steady state, or grow in amplitude leading to thermal oscillations or uncontrollable thermal runaway. When stability analysis is carried out over a range of the design and operational parameters of the system it is often called bifurcation analysis, and it can provide a valuable stability map of the system.

The mathematical theory and methodology behind stability analysis are well-established and accepted. Briefly, at each steady state solution of the parent dynamical system we construct a dynamical system for  the perturbation to that steady state. Since the perturbation is small we can write it as a Taylor series expansion and retain only the first-order terms. According to the theory
of linear differential equations, the solution can be written as a superposition of
terms of the form $e^{\lambda_j}t$ where $\{\lambda_j\}$ is the set of eigenvalues of the Jacobian matrix of coefficients of the linearized perturbation system. A nonzero complex part of an eigenvalue contributes an oscillatory component to the solution. If the real part of an eigenvalue is positive the perturbation must grow exponentially  with time. A stable steady state, therefore, is one for which \textit{all} the eigenvalues of the Jacobian  of the linearized perturbation system have negative real components. 

\setlength\parskip{0mm}
\section{Methodology and dynamical model\label{sec3}}
\subsection{Methodology}
The following general procedure for stability analysis of thermoreactive systems was applied: 
\begin{enumerate}
\item Write down dynamical coupled mass and enthalpy balances for the reaction system.
\item Choose a primary \textit{bifurcation parameter}, usually a parameter that can be tuned experimentally, such as an inlet temperature or flow rate.
\item Compute the steady state solutions as a function of the bifurcation parameter.
\item Evaluate the stability of each solution by computing the characteristic  eigenvalues of the linear perturbation at each point.
\item Flag each change in sign of the real parts of the eigenvalues. These are the singular points. 
\item Compute the amplitude, period, and stability of any periodic solutions as a function of the bifurcation parameter. 
\item If appropriate, repeat the analysis using another bifurcation parameter. If singular points are found  compute their loci using a second parameter to obtain a stability map of the dynamical system. 
\end{enumerate}
This modelling and analysis was carried out in the well-stirred approximation. Results from this analysis  provide essential guidance for design and operation of an economically and safety optimised reactor system, and may be used to inform expensive convective
simulations that require substantial high-performance computational resources.

\setlength\parskip{2mm}

\subsection{Endex-coupled carboniser-calciner: dynamical model}

Enthalpy summation and mass balances for the gas-phase reactant lead to the following dynamical equations for the Endex carboniser-calciner system: 
\begin{align}
V_1\frac{dc_1}{dt} = & -V_1 v_1(T_1,p_1)
 + F_1(c_{1,\rm{in}} - c_1)\label{e1}\\[1mm]
V_1\overline{C}_1\frac{dT_1}{dt} =  &\: V_1(-\Delta H)v_1(T_1,p_1) 
  + F_1 \overline{C}_{1,\rm{g}}(T_{1,\rm{in}} - T_1)
				+ \left(F_sC_s+L_{\rm{ex}}\right)\left(T_2 - T{_1}\right)
				 \label{e2}\\[2mm]
V_2\frac{dc_2}{dt} = & V_2 v_2(T_2,p_2) - F_2c_2 \label{e3}\\[1mm]
V_2\overline{C}_2\frac{dT_2}{dt} = & \: V_2\Delta H v_2(T_2,p_2)
 - F_2 \overline{C}_{2,\rm{g}}T_2
				+ \left(F_sC_s+L_{\rm{ex}}\right)\left(T_1 - T_2\right).\label{e4}
				\end{align}
 The quantities and notation are defined in table \ref{table1}, Appendix A. Equations (\ref{e1})--(\ref{e4}) describe an Endex system that is heat-coupled through cycling of loaded and unloaded sorbent between the carboniser (denoted by subscript~1) and calciner (denoted by subscript 2), and via direct common wall transfer. The overall system is treated as fully insulated. In reality it is expected that heat losses to the environment will be small, less than around 5\%. 
 
The solid-gas surface reaction rates are functions of temperature $T_i$, pressure $p_i$, and  fractional surface coverage $\theta_i$, $(i=1,2)$, 
\begin{align} 
v_1(T_1,p_1)=&\left(\tfrac{p_1}{p_{1,\rm{eq}}} - 1 \right)\theta_1 \epsilon k(T_1) \, \zeta_1 \,S \label{e5}\\
v_2(T_2,p_2)=&\left(1-\tfrac{p_2}{p_{2,\rm{eq}}}  \right)
\left(1-\theta_2\right) \epsilon k(T_2) \, \zeta_2 \,S ,\label{e6}
\end{align}
where the rate constant $k(T_i)$ has the usual Arrhenius temperature dependence, and
$\theta_i$ is given by the Langmuir isotherm
\begin{equation}
\theta_i = 
\tfrac{\left(\tfrac{p_i}{p_{i,\rm{eq}}}\right)^{1/2}}
{1 + \left(\tfrac{p_i}{p_{i,\rm{eq}}}\right)^{1/2}}. \label{e7}
\end{equation} 
In this model, the isotherm is based on the CO$_2$ molecule occupying two surface sites and the saturation pressure of CO$_2$ is pinned to the equilibrium pressure of CO$_2$ in the sorbent.        

 \subsection{Data source}
Physicochemical data and reaction rates for the surface reaction system CaO$_{\rm(s)}$ + CO$_{2\rm(g)}$ $\rightleftarrows$ CaCO$_{3\rm(s)}$ were taken from the published literature and tabulated in \cite{Sceats:2009b}.   The reactor design parameters and solids and gas flow rates are those given in \cite{Sceats:2009b} for a demonstration unit for scrubbing the emissions from a 5\,MW lignite fuelled power plant. Numerical values of data and quantities used in this analysis are given in table \ref{table1}, Appendix A.

\section{Results and discussion\label{sec4}}

\subsection{Standalone carboniser}
\begin{figure}[t]
\centerline{
  \includegraphics[scale=0.56]{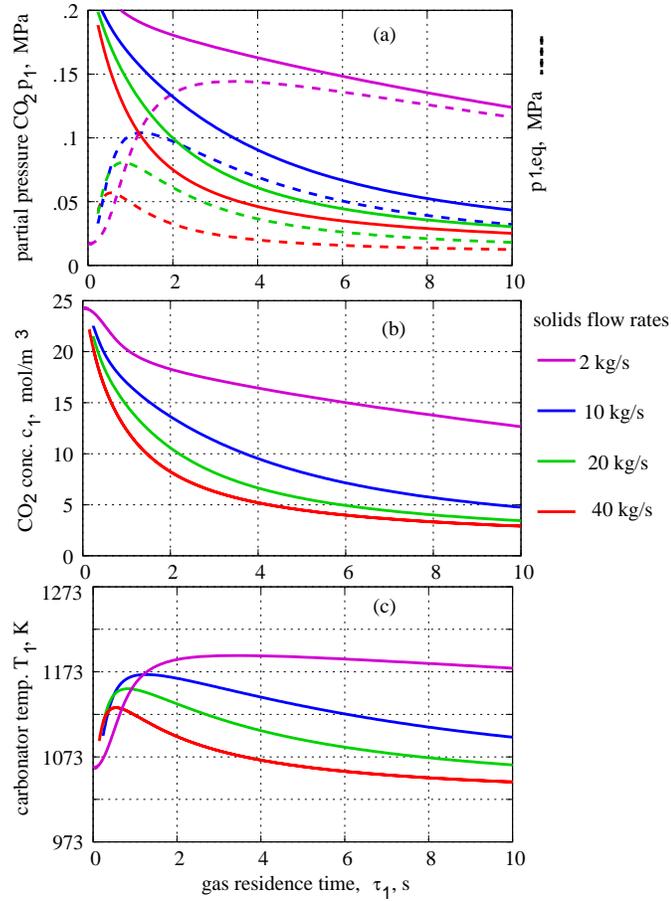} }
  \caption[]{\label{fig2} Steady state analysis of the standalone carboniser. $T_{1,\rm{in}}= 1060$\,K, $T_{s,\rm{in}}= 1021$\,K. }
\end{figure}

The first step in the analysis of the Endex dynamical system, equations (\ref{e1})--(\ref{e4}), is to analyse the  carboniser segment in standalone mode. The dynamics of this system will be a subset of the Endex dynamics. The standalone carboniser is described by the following smaller, simpler dynamical system:
\begin{align}
V_1\frac{dc_1}{dt} = & -V_1 v_1(T_1,p_1)
 + F_1(c_{1,\rm{in}} - c_1)\tag{\ref{e1}}\\[1mm]
V_1\overline{C}_1\frac{dT_1}{dt} =  &\: V_1(-\Delta H)v_1(T_1,p_1) 
  + F_1 \overline{C}_{1,\rm{g}}(T_{1,\rm{in}} - T_1)
				+ F_sC_s\left(T_{s,\rm{in}} - T{_1}\right). 
				 \label{e8}
\end{align}	
Steady state solutions of equations (\ref{e1}) and (\ref{e8}), using equation \ref{e5}, were computed as a function of the gas residence time, $\tau_1\equiv V_1/F_1$, and are rendered in figure \ref{fig2} in terms of the dynamical variables $p_1$ (a), $c_1$ (b), and $T_1$ (c). The computations were carried out for four values of the solids flow rate $F_s$ as indicated. The real parts of the two eigenvalues of the linear perturbation were positive at each point, thus the steady states are stable within this regime.
 
 In (a) the equilibrium partial pressure, $p_{1,\rm{eq}}$, corresponding to each $F_s$ has also been plotted as dotted lines. We see that the condition $p_1\gg p_{1,\rm{eq}}$ holds over the range of residence times considered, thus carbonation proceeds spontaneously. 
 
  In (b) we see that higher mass flow rates of solid sorbent allow improved uptake of the CO$_2$ by the sorbent at shorter gas residence times. For example, to achieve a CO$_2$ concentration of $\sim$7\,mol/m$^3$ (corresponding to a CO$_2$ conversion of 71\%) requires a residence time of $\sim$7.2\,s at a solids flow rate of 10\,kg/s, but only 4\,s if the solids flow rate is 20\,kg/s. 
 
However, if it is desired to maintain the temperature below  $\sim$1123\.K and achieve appreciable CO$_2$ conversion we see from (c) that gas residence times~$>$~6\,s and solids flow rates $>$ 10\,kg/s are required.

\clearpage

\subsection{Endex-coupled carboniser-calciner: analysis}
In the standalone carboniser the heat carried by the partially carbonated sorbent was discarded into the environment. In effect the standalone carboniser loses the heat generated by the adsorption reaction to a heat bath held at the constant temperature $T_{s,\rm{in}}$ at the rate indicated by the last term in equation (\ref{e8}). 

In the Endex-coupled system, modelled by equations (\ref{e1})--(\ref{e4}), the heat generated by adsorption and carried by the partially carbonated sorbent is recovered directly to drive the calcination of the sorbent. There is no external thermostat or heat bath; instead, the dynamical temperature of the carboniser $T_1$ and that of the calciner $T_2$ are coupled at the rate indicated by the last term in equations~(\ref{e2}) and (\ref{e4}).

\subsubsection*{Zero wall heat exchange}
 Assuming  in the first instance that $L_{\rm{ex}}=0$, i.e., the carboniser and calciner communicate thermally only via transfer of the sorbent, we have the CO$_2$ residence time in the calciner $\tau_2\equiv V_2/F_2$ as the only additional tunable parameter.  

\begin{figure}[hb]
\centerline{
  \includegraphics[scale=0.56]{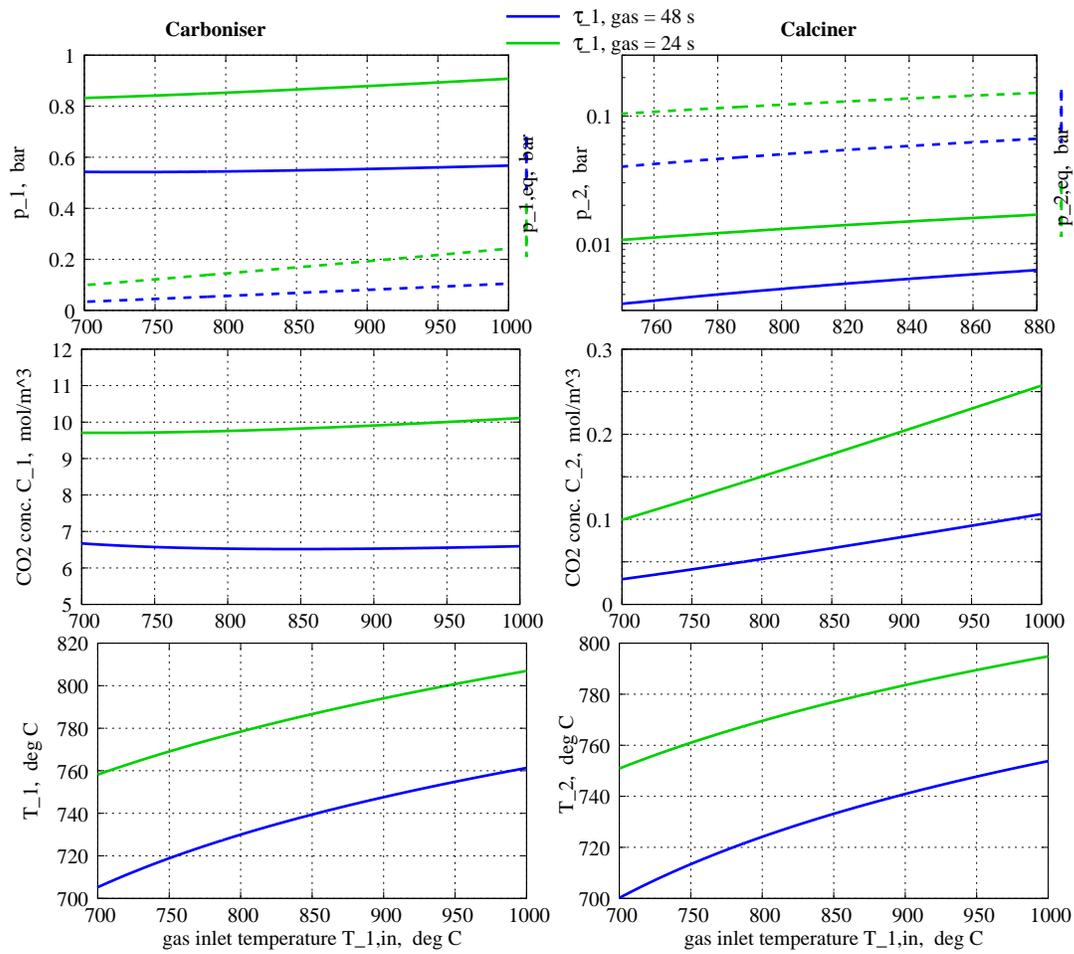}}
  \caption[]{\label{fig3}Left subfigures: carboniser, right subfigures: calciner.  $\tau_{2,\rm{gas}}=$ 30\,s.}
\end{figure}

\begin{figure}[ht]
\centerline{
  \includegraphics[scale=0.56]{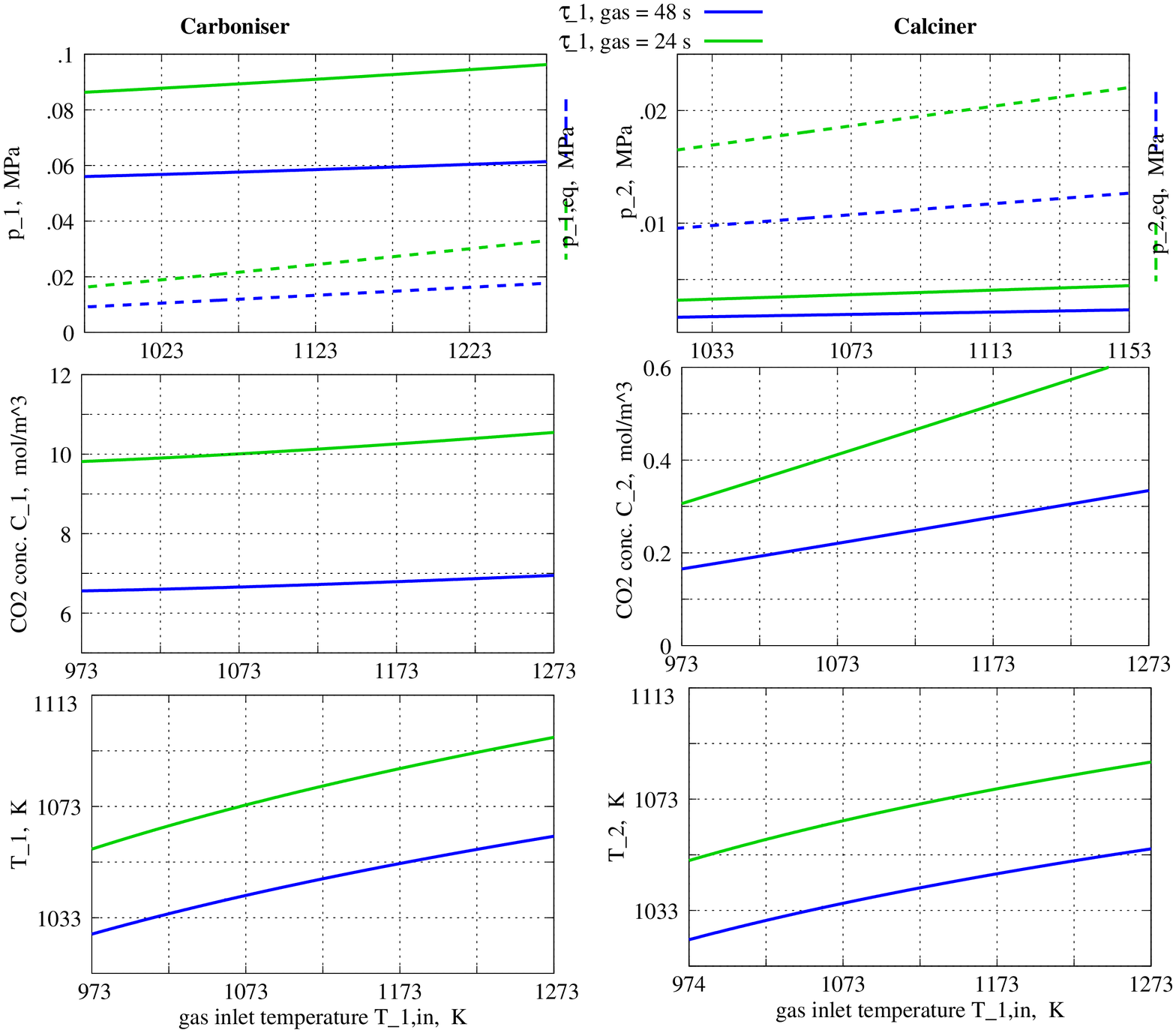}}
  \caption[]{\label{fig4}Left subfigures: carboniser, right subfigures: calciner. $\tau_{2,\rm{gas}}=$ 60\,s.}
\end{figure}

The first task is to select an optimum range for the gas inlet temperature, $T_{1,\rm{in}}$. The putative 
set-point for $T_{1,\rm{in}}$ is around 1023\,K. In figures \ref{fig3} and \ref{fig4} steady state solutions of equations \ref{e1}--\ref{e4} are plotted with $T_{1,\rm{in}}$ as the bifurcation parameter. The eigenvalue analysis gave the steady states as stable over this range.

In figure \ref{fig3} the calciner gas residence time $\tau_{2,\rm{gas}}$ is set at 30\,s, and and in figure \ref{fig4} $\tau_{2,\rm{gas}}$ is set at 60\,s. The steady state solutions have been computed and plotted for two values of the carboniser gas residence time $\tau_{1,\rm{gas}}$. 

It is interesting to observe that the longer carboniser gas residence time $\tau_{1,\rm{gas}}$ (blue lines) gives much improved CO$_2$ uptake by the sorbent over the shorter  $\tau_{1,\rm{gas}}$ (green lines), yet the dynamical carboniser temperature $T_1$ is depressed. One normally expects that an exothermic reacting system will become hotter for higher conversion of the reactant. The occurrence of the  back  reaction (calcination) at higher conversion of the CO$_2$ should work against this, but the positive enthalpy change associated with the back reaction is not built into equation \ref{e4}. In this case it is the high rate of thermal transport provided by the mass flow of sorbent with high thermal capacitance that allows this ``temperature inversion''. When the sorbent mass flow rate $F_s$ is set to zero, the carboniser and calciner are decoupled and the normal increase in $T_1$ with longer gas residence time occurs. 
\clearpage

\begin{figure}[ht]
\centerline{
  \includegraphics[scale=0.6]{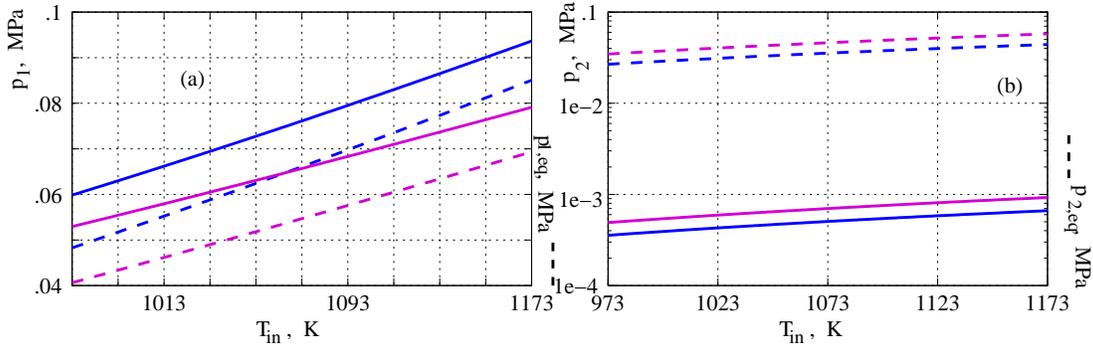}}
  \caption[]{\label{fig5} Blue: $F_s= 10$\,kg/s, magenta: $F_s= 40$\,kg/s. $\tau_{1,\rm{gas}}=$ 15\,s, $\tau_{2,\rm{gas}}=$ 15\,s, $L_{\rm{ex}}=0$. }
\end{figure}

What this means is that the sorbent mass flow is an important control tool for the system. We now know, for example, that $\tau_{1,\rm{gas}}$  can be made as long as you like without risking overheating of the carboniser, provided $F_s$ is maintained above some critical rate. In figure \ref{fig5} the partial pressure steady states are compared for two values of $F_s$. The higher sorbent flow rate evidently gives better performance. 

Comparing figures \ref{fig3} and \ref{fig4} we note another interesting effect of thermal coupling of the reactors: the carboniser ``sees'' the calciner gas residence time  $\tau_{1,\rm{gas}}$. Since the negative enthalpy change of the back reaction (carbonation) is not built into equation \ref{e4}, it is the longer $\tau_{1,\rm{gas}}$ that results in warmer sorbent entering the carboniser at $T_2$, which in turn increases $T_1$. This improved thermal coupling decreases $T_1-T_2$, consistent with the rigorous ideal Endex result \citep{Gray:1999}: $\lim_{\sigma\rightarrow 0} |T_1-T_2| = 0$, where $\sigma=1/F_s$.  

For both figures, in the carboniser segment the condition $p_1\gg p_{1,\rm{eq}}$ holds and in  
the calciner segment the condition  $p_2\ll p_{2,\rm{eq}}$ holds. The solutions are stable over the computed regime. However, dynamical stability---i.e., the behaviour of complex conjugate pairs of eigenvalues, or the occurrence of Hopf bifurcations---is governed by the volumetric specific heats. A system with low thermal capacitance is likely to become thermally unstable, and in fact Hopf bifurcations do occur in this system if  artificially low specific heats are used. 
		 
\begin{figure}[t]
\centerline{
  \includegraphics[scale=0.75]{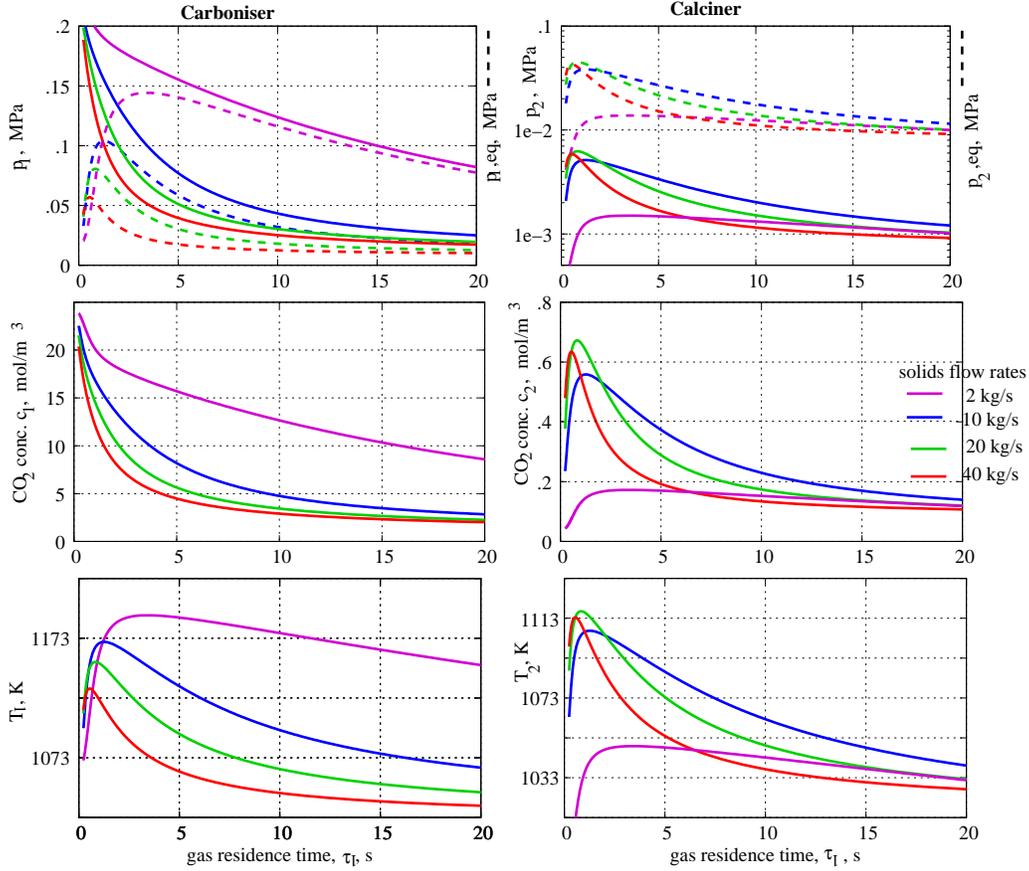}}
  \caption[]{\label{fig6}Left subfigures: carboniser, right subfigures: calciner. Endex steady states as a function of carboniser gas residence time $\tau_{1,\rm{gas}}$. $\tau_{2,\rm{gas}}=$ 30\,s, $L_{\rm{ex}}=0$.}
\end{figure}

The steady state solutions with $\tau_1$ as bifurcation parameter are plotted in figure \ref{fig6}, for four solids flow rates. In the carboniser a residence time greater than 10\,s is required to achieve a CO$_2$ uptake of 90\% or more. However, residence times longer than about 15\,s do not lead to appreciably more uptake of CO$_2$.

\subsubsection*{Non-zero wall heat transfer}
\begin{figure}[t]
\centerline{
  \includegraphics[scale=0.6]{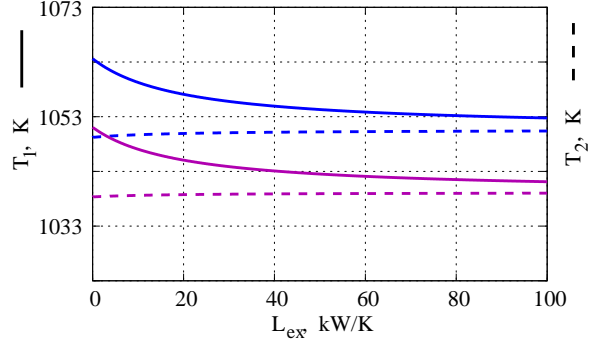}}
  \caption[]{\label{fig7} Blue lines: $\tau_{1,\rm{gas}}=$ 10\,s; Magenta lines: $\tau_{1,\rm{gas}}=$ 15\,s. $T_1$ is plotted with solid lines, $T_2$ with dashed lines. $\tau_{2,\rm{gas}}=$ 30\,s, $F_s=20$\,kg/s.}
\end{figure}
\begin{figure}
\centerline{
  \includegraphics[scale=.6]{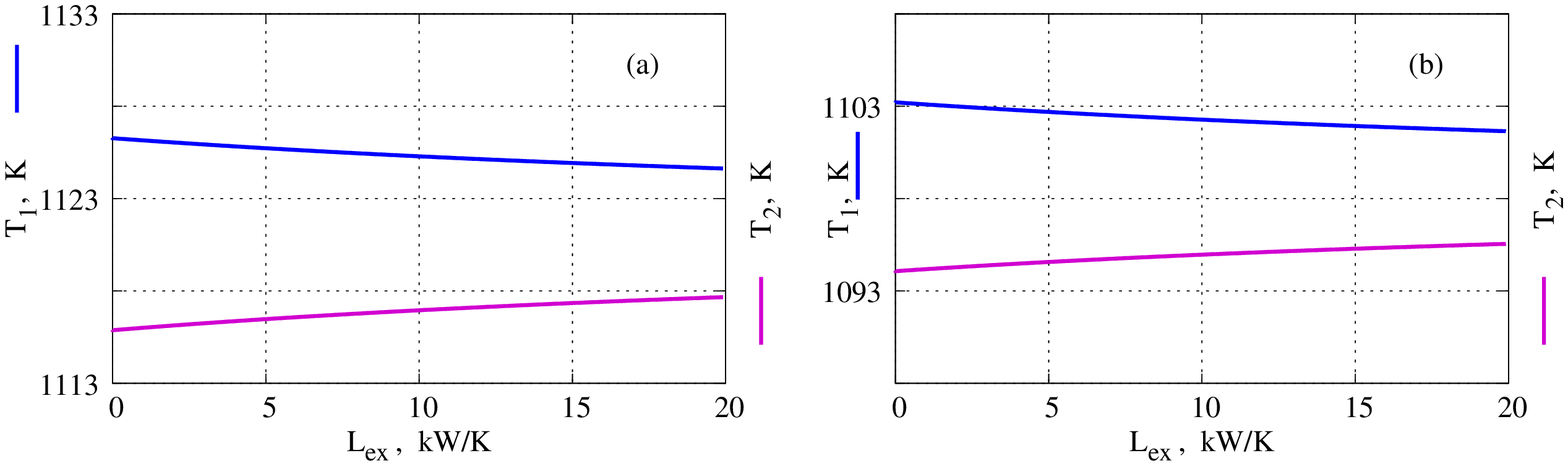}}
  \caption[]{\label{fig8} (a)  $\tau_{1,\rm{gas}}=$ 15\,s, (b) $\tau_{1,\rm{gas}}=$ 20\,s. $\tau_{2,\rm{gas}}=$ 30\,s, $F_s=40$\,kg/s. $T_1$ is plotted with blue lines, $T_2$ with magenta lines.}
\end{figure}

In this scenario the carboniser and calciner share a common wall or walls, which are also heat exchangers. For cylindrical reactor segments the carboniser cylinder would be embedded in the calciner segment. (The volume displacement is negligible and has not been corrected for in this analysis. 

For the purpose of simulating the effects of wall heat exchange we now regard the combined heat transfer coefficient $L_{\rm{ex}} $ as a tunable bifurcation parameter, and the steady state solutions are plotted in figures \ref{fig7} and \ref{fig8}. The theoretical limit as $1/L_{\rm{ex}} \rightarrow 0$ is $T_1=T_2$. However a reasonable upper limit of $L_{\rm{ex}} $ for this configuration would be around 10\,kW/K. It will be shown in subsection \textbf{C} below that nonzero wall heat transfer is desirable from safety considerations.

\subsubsection*{ Start-up and shut-down dynamics}
\begin{figure}[hb]
\centerline{
  \includegraphics[scale=.6]{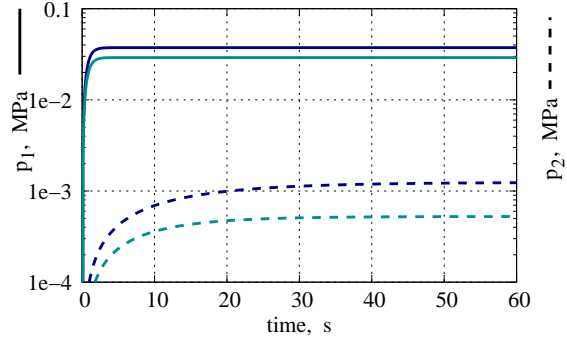}}
  \caption[]{\label{fig9} Cyan coloured lines:  $\tau_{2,\rm{gas}}=$ 10\,s. Indigo coloured lines: $\tau_{2,\rm{gas}}=$ 10\,s. In both cases $\tau_{1,\rm{gas}}=$ 15\,s, $F_s=20$\,kg/s, $L_{\rm ex}=0$. In each case the initial CO$_2$ concentrations are zero and the initial temperatures are set equal to the steady state temperatures.}
\end{figure}

\textit{Start-up:} The trajectories plotted in figure \ref{fig9} indicate that full steady state operation can be achieved in less than 60\,s, whiler quasi steady-state operation of the carboniser segment is achieved in less than 5\,s. The time lag is due to the much higher effective activation energy for the calcination, and it could be reduced by employing operational strategies that achieve improved effective kinetic matching of the reactions, for example, increasing the pumpout rate~$F_2$. 

\begin{figure}[ht]
\centerline{
  \includegraphics[scale=0.7]{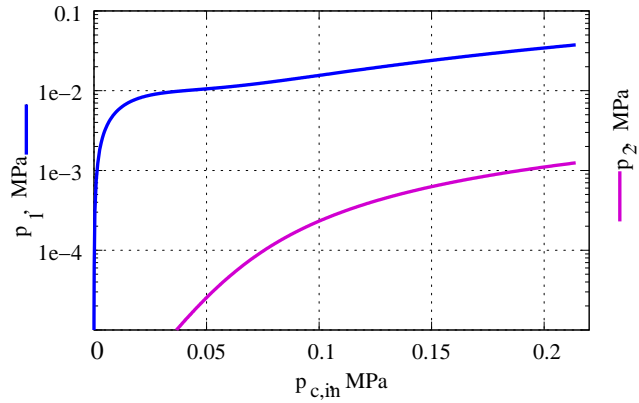}}
  \caption[]{\label{fig10} The gas inlet pressure is allowed to ebb quasistatically. $F_s=20$\,kg/s,  $\tau_{1,\rm{gas}}=$ 15\,s, $\tau_{2,\rm{gas}}=$ 15\,s, $T_{1,\rm{in}} = 1060$\,K, $L_{\rm{ex}}=0$. }
\end{figure}

\textit{Shutdown:}  In a normal shutdown scenario the carboniser-calciner shutdown dynamics are coupled to the shutdown dynamics of the flue gas generator. Assuming the CO$_2$ partial pressure $p_{\rm{c,in}}$ at the inlet is reduced linearly and quasistatically we may simulate this type of shutdown as a quasi steady-state procedure,  using $p_{\rm{c,in}}$  as the bifurcation parameter. In figure \ref{fig10} the carboniser partial pressure $p_1$ ebbs slowly until $p_{\rm{c,in}}$ is about 0.025\,MPa then drops off dramatically. The calciner gas pressure $p_2$ declines smoothly as the loading of the sorbent declines. 

\clearpage

\begin{figure}[h]
\centerline{
  \includegraphics[scale=0.7]{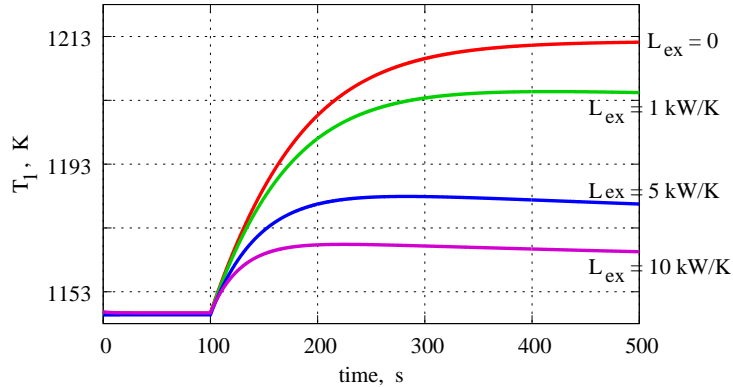}}
  \caption[]{\label{fig11} At $t=100$\,s the solids flow rate $F_s$ is switched from 40\,kg/s to zero.  $\tau_{1,\rm{gas}}=$ 15\,s, $\tau_{2,\rm{gas}}=$ 15\,s, $T_{1,\rm{in}} = 1060$\,K. }
\end{figure}

\textit{Possible danger:} For an Endex carboniser-calciner reactor system without wall thermal contact the only channel for heat transfer to the calciner is via the sorbent. If the sorbent flow $F_s$ is interrrupted (perhaps by mechanical failure) the cooling capacity of the calciner is unavailable. The only channel for heat removal from the carboniser is then via the outflow of scrubbed flue gas. The adiabatic temperature rise for complete conversion of 
24.28\,mol/m$^{-3}$ of CO$_2$ to carbonate is around 650\,K. 
Although a temperature rise of this magnitude would not in practice occur, because the scrubbed gas outflow would continue to remove heat, the scenario in which $F_s$ is interrrupted is of concern. It is simulated in figure \ref{fig11}, in which steady state operation is interrupted by switching off the solids flow.  The time series were computed for zero wall heat transfer, $L_{\rm{ex}}=0$, and for $L_{\rm{ex}}=1$, 5 and 10 kW/K. The thermal excursion is large and dramatic for $L_{\rm{ex}}=0$, more than 80\,K and possibly exceeding the safety limits of the vessel. Where the calciner and carboniser have efficient wall thermal contact, however, the thermal excursion is much smaller and quite manageable. In the case of $L_{\rm{ex}}=$ 10 kW/K the maximum temperature increase of 15\,K occurs 100\,s after the solids flow is switched off, after which the temperature declines slowly. 

Thus the inclusion of efficient wall thermal contact in the Endex reactor design may be an important safety consideration. 

\subsection*{D\quad Instabilities in the system}
In the subsections above we found that within a broad range of the expected normal operationg conditions the Endex carboniser-calciner reactor is free of instabilities. This raises the question of where, exactly, in the parameter space may instability occur in this system, since there are no global constraints on the nonlinear dynamical system, equations (\ref{e1})--(\ref{e4}) that forbid the occurrence of a positive real eigenvalue part. 

In fact multiplicity and hysteresis may occur at low temperatures of the gas inlet.
\begin{figure}[ht]
\centerline{
  \includegraphics[scale=0.7]{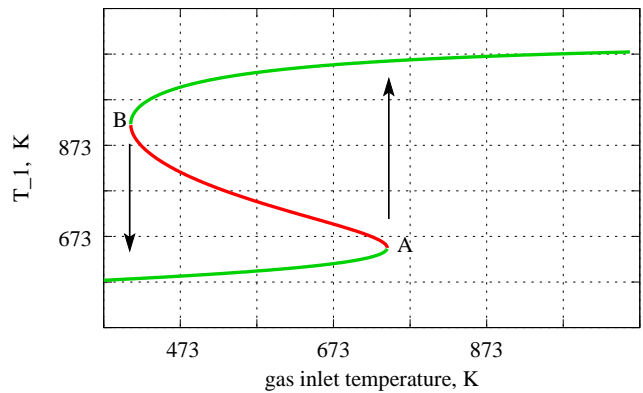}}
  \caption[]{\label{fig12} Green curves are the stable steady states, red curves are the unstable steady states. $F_s=5$\,kg/s,  $\tau_{1,\rm{gas}}=$ 2.4\,s, $\tau_{2,\rm{gas}}=$ 15\,s, $L_{\rm{ex}}=0$. }
\end{figure}
  In figure \ref{fig12} the steady states have been computed and plotted as a function of the gas inlet temperature $T_{1,\rm{in}}$. The eigenvalues were also monitored, and the stability of the steady states is colour-coded. 
  
One may, in principle, carry out an experiment in which the gas inlet temperature is increased quasistatically from, say, 473\,K. On this lower green curve the reaction temperature $T_1$ increases slowly but small perturbations decay. When the gas inlet temperature reaches the value marked A on the figure the solutions change character. A real eigenvalue passes through  zero and becomes positive. Small perturbations now grow in amplitude. A rapid transition to the high-temperature stable branch occurs, indicated by an arrow on the figure. The inlet temperature may then be tuned up to its normal set-point at 1060\,K, or it may be tuned down along the upper green curve.  When it reaches the point marked B on the figure an abrupt temperature collapse must occur, since the solutions again lose stability. 

Although this hysteresis loop may be exploited in start-up and shut-down procedures, it can be seen that it is not relevant over the normal operating regime of the reactor. 

Thermal oscillations associated with nonzero complex components of the eigenvalues do not occur in this system. This is because such behaviour is typically governed by the dynamics of the exothermic reaction, i.e., the left hand side of equation (\ref{e2}).  The high thermal capacitance provided by the solids fraction effectively damps any ocillatory components, or complex parts of the eigenvalues.

\section{Summary and conclusions\label{sec5}}

A dynamical system model was derived for an Endex coupled carboniser-calciner housing the CaO/CaCO$_3$ surface reactions, in the well-stirred fully insulated approximation. 

Steady state and stability analysis of the carboniser compartment in standalone mode provided a subset of states from which to begin the more difficult task of analysing the full Endex system. Approximate lower bounds for the gas residence time $\tau_1$ and solids flow rate $F_s$ were obtained. 
In the steady state and stability analysis of the full Endex system it was observed that the sorbent mass flow rate is an important control parameter, because the thermal transport from carboniser to calciner provided by the sorbent depresses the carboniser temperature at long gas residence~times. 

The equilibrium partial pressures in carbonsiser and calciner segments were monitored. Over the parameter regime studied $p_{1,\rm{eq}}<p_1$ and $p_{2,\rm{eq}}>p_2$, as required by the surface reaction thermokinetics for spontaneous carbonation and calcination respectively.  

 The parameter regime studied is given in the following table:\\

\centerline{
\begin{tabular}{p{0.1\textwidth}p{0.2\textwidth}}
\hline
$T_{1,\rm{in}}$ & 973\,--\,1273\,K\\
$\tau_1$ & 0.1\,--\,20\,s\\
$F_s$& 10\,--\,40\,kg/s\\
$\tau_2$ & 15\,--\,60\,s\\
$L_{\rm{ex}}$ & 0\,--\,100\,kW/K\\
\hline\\
\end{tabular}
}

The real parts of the eigenvalues of the solutions remain negative over the parameter regime studied, hence steady state operation is stable over these ranges; i.e., perturbations to the system set-point decay rather than grow.

Start-up time to steady state operation was found to be less than 60\,s for initial CO$_2$ partial pressure of zero and initial temperatures equal to the steady state temperatures. Quasi steady state operation of the carboniser segment can be achieved in less than 5\,s. The shut-down dynamics of the reactor were modelled as the gradual decline in pressure of the gas inflow. The carboniser partial pressure fell slowly until the temperature became too low to sustain the reaction, whereupon the partial pressure drops off rapidly.

Interruption of the sorbent flow was flagged as a possible source of dynamical thermal instability, in a reactor with no or poor wall thermal coupling.  It was found that a reactor design with significant wall heat exchange between the carboniser and calciner compartments insured against temperature surges in the carboniser in the event of an interruption to the solids flow.

To complete the stability analysis, the location of instabilities in this Endex system was pinpointed as a region of thermal multiplicity over temperatures well below the normal operating regime. It was noted that oscillatory instabilities are fully damped by the high thermal capacitance provided by the sorbent. 

These modelling and analysis results confirm that the proposed reactor configuration for the Calcium Looping reactor is a stable Endex configuration that can, in principle, scrub CO$_2$ from a flue gas stream efficiently and regenerate a pure stream of CO$_2$ for geosequestration without additional energy requirements. The system has the potential to capture more than 90\% of the CO$_2 $ from flue gas emissions and release it in a pure stream, in a thermally safe reactor that requires no thermostatting.  The system is stable to perturbations and exhibits gas switching during start-up.  

\clearpage

\appendix
\section{Notation}
\setlongtables
{\footnotesize
\begin{longtable}{p{1cm}p{8cm}l}\caption{For compactness the subscript $i$ is used as appropriate, where $i=1,2$; 1 refers to the carboniser and 2 refers to the calciner. \label{table1} }
\\\hline\hline\endhead
&&\\
$c_i$ & concentration of CO$_2$ & mol/m$^3$\\
$c_{1,\rm{in}}$ & inflow concentration of CO$_2$& $p_{\rm{c,in}}/RT_{1,\rm{in}}=24.3$\,mol m$^{-3}$\\
$k(T_i)$ & rate constant & $114\exp(-E/RT_1) \rm{mol}/(\rm{m}^3 \rm{s}^1)$ \\
$p_0 $& & 4.147e06\,MPa\\
$p_{\rm{c,in}}$& inlet partial pressure CO$_2$ &MPa \\
$p_i$ & partial pressure CO$_2$ &MPa \\
$p_{i,\rm{eq}}$& equilibrium partial pressure of CO$_2$& $p_0 \exp(-|\Delta H|/RT_i)$  MPa\\
$v_i$ & reaction rate & mol\,m$^{-3}$s$^{-1}$\\
$\overline{C}_1$& weighted volumetric specific heat of carboniser contents & 160\,$\text{kJ K}^{-1}\text{m}^{-3}$\\
$\overline{C}_{1,g}$& weighted volumetric specific heat of gas& 5.8 $\text{kJ K}^{-1}\text{m}^{-3}$\\
$\overline{C}_2$& weighted volumetric specific heat of calciner contents & $25 \,\text{kJ K}^{-1}\text{m}^{-3}$\\
$\overline{C}_{2,g}$& volumetric specific heat of calciner gas& 25 $\text{J K}^{-1}\text{m}^{-3}$\\
$C_s$& specific heat of sorbent& $975   \text{J\,K}^{-1}\text{kg}^{-1}$\\
$E$ & activation energy for calcination & 205 kJ/mol\\
$F_1$ & volumetric flow rate of gas into the carboniser & $V_1/\tau_1$ $\text{m}^3\text{s}^{-1}$\\
$F_2$ & volumetric outflow rate of gas from the calciner & $V_2/\tau_2$ $\text{m}^3\text{s}^{-1}$\\
$F_s$ & mass flow rate of sorbent & $\text{kg s}^{-1}$ \\
$L_{\rm{ex}}$ & heat exchange rate coefficient&  kW\,K$^{-1}$\\
$R$& gas constant & 8.314 J/(mol K)\\
$S$ & surface area & $5e07\, \text{m}^2/\text{m}^3$ \\
$T_i$ & reactor segment temperature  & K \\
$T_{1,\rm{in}}$& temperature of gas at the inlet& K\\
$V_1$ & internal volume of carboniser & $\pi\times 0.25^2\times 12= 2.356\, \text{m}^3$\\
$V_2$ & internal volume of  calciner& $\pi\times 2^2\times 12= 150.8\, \text{m}^3$\\
$\Delta H$& reaction enthalpy & $-170 $\,kJ/mol CO$_2$\\
$\epsilon$ & porosity of nascent lime & 0.51\\
$\tau_i$ & gas residence time&s\\
$\theta_i$ & fractional sorbent surface coverage&\\
$\zeta_1 $ & carboniser solid fraction & 0.5 \\
$\zeta_2 $ & calciner solid fraction & 0.008 \\
\hline\hline

\end{longtable}}

\noindent\textit{Acknowledgement:} This work is supported by Australian Research Council Future Fellowship FT0991007 (R. Ball). We thank an anonymous referee for helpful and constructive comments.
\clearpage

\bibliographystyle{elsarticle-num}

\end{document}